\documentclass{gtart_h}

\def\ifplaintex{\expandafter\ifx\csname documentclass\endcsname\relax}

\def\gtp{{\mathsurround=0pt\it $\cal G\mskip-2mu$eometry \&\ 
$\cal T\!\!$opology $\cal P\!$ublications}}  

\def\recd{{\small Received:\qua\receiveddate\ifx\reviseddate\relax
\else\qquad Revised:\qua\reviseddate\fi\par}} 

\def\lognumber#1{\def\thelognumber{#1}}
\def\volumenumber#1{\def\thevolumenumber{#1}}
\def\volumeyear#1{\def\thevolumeyear{#1}}
\def\papernumber#1{\def\thepapernumber{#1}}
\def\pagenumbers#1#2{\def\startpage{#1}\def\finishpage{#2}}
\def\published#1{\def\publishdate{#1}}

\def\received#1{\def\receiveddate{#1}}
\def\revised#1{\def\reviseddate{#1}}
\def\accepted#1{\def\accepteddate{#1}}

\def\asciiauthors#1{\def\theasciiauthors{#1}}
\def\asciiaddress#1{\def\theasciiaddress{#1}}

\def\coverauthors#1{\def\thecoverauthors{#1}}

\let\\\par\let\thelognumber\relax\let\thevolumenumber\relax
\let\thepapernumber\relax\let\thevolumeyear\relax\let\startpage\relax
\let\finishpage\relax\let\publishdate\relax\let\receiveddate\relax
\let\reviseddate\relax\let\accepteddate\relax\let\theasciititle\relax
\let\theasciiauthors\relax\let\theasciiaddress\relax
\let\theasciiabstract\relax

\let\thecoverauthors\relax\let\theasciiemail\relax

\ifplaintex
\font\logobig=cmssbx10 scaled 3836
\font\logomed=cmssbx10 scaled 2557
\else
\font\logobig=cmssbx10 scaled 4200
\font\logomed=cmssbx10 scaled 2800
\fi

\long\def\makeagttitle{   
\count0=\startpage
\agt\hfill      
\hbox to 45truept{\vbox to 0pt{\vglue -13truept{\logomed A\kern -.37em{\logobig 
T}\kern -.38em G}\vss}\hss}
\break
{\small Volume \thevolumenumber\ (\thevolumeyear)
\startpage--\finishpage\nl
Published: \publishdate}

\vglue .25truein

{\parskip=0pt\leftskip 0pt plus
1fil\def\\{\par\smallskip}{\Large\bf\thetitle}\par\medskip} \vglue
0.05truein

{\parskip=0pt\leftskip 0pt plus 1fil\def\\{\par}{\sc\theauthors}
\par\medskip}%
 
\vglue 0.03truein 

{\small\leftskip 25truept\rightskip 25truept{\bf Abstract}\stdspace\theabstract

{\bf AMS Classification}\stdspace\theprimaryclass
\ifx\thesecondaryclass\relax\else; \thesecondaryclass\fi\par
{\bf Keywords}\stdspace \thekeywords\par}\vglue 7truept

}   

\ifplaintex
\font\phead=cmsl9 scaled 950
\font\pnum=cmbx10 scaled 913
\font\pfoot=cmsl9 scaled 950
\headline{\vbox to 0pt{\vskip -4.5mm\line{\small\phead\ifnum
\count0=\startpage ISSN 1472-2739 (on-line) 1472-2747 (printed)
\hfill {\pnum\folio}\else\ifodd\count0\def\\{ }%
\ifx\theshorttitle\relax\thetitle\else\theshorttitle\fi\hfill{\pnum\folio}
\else\def\\{ and }{\pnum\folio}\hfill\ifx\theshortauthors\relax\theauthors
\else\theshortauthors\fi\fi\fi}\vss}}
\footline{\vbox to 0pt{\vglue 0mm\line{\small\pfoot\ifnum\count0=\startpage
\copyright\ \gtp\hfill\else
\agt, Volume \thevolumenumber\ (\thevolumeyear)\hfill\fi}\vss}}
\else
\font=cmsl9 scaled 1050
\font\lnum=cmbx10 
\font=cmsl9 scaled 1050
\makeatletter
\def\@oddhead{{\small\lhead\ifnum\count0=\startpage ISSN 1472-2739 
(on-line) 1472-2747 (printed)\hfill {\lnum\number\count0}\else\ifodd\count0
\def\\{ }\ifx\theshorttitle\relax \thetitle \else\theshorttitle\fi\hfill
{\lnum\number\count0}\else\def\\{ and }{\lnum\number\count0}
\hfill\ifx\theshortauthors\relax 
\theauthors\else\theshortauthors\fi\fi\fi}}\def\@evenhead{\@oddhead}
\def\@oddfoot{\small\lfoot\ifnum\count0=\startpage\copyright\ \gtp\hfill\else
\agt, Volume \thevolumenumber\ (\thevolumeyear)\hfill\fi}
\def\@evenfoot{\@oddfoot}
\makeatother
\fi
\let\maketitlepage\makeagttitle
\let\makeshorttitle\maketitlepage
\let\maketitle\maketitlepage

\newwrite\gtoutfile
\long\gdef\makeheadfile{  
{\def\\{, }\def\s{ }
\immediate\openout\gtoutfile head.xxx
\immediate\write\gtoutfile{Proxy-for: \ifx\theasciiauthors\relax
\theauthors\else\theasciiauthors\fi\s<\ifx\theasciiemail\relax\theemail\else\theasciiemail\fi>}
\immediate\write\gtoutfile{\noexpand\\}
\immediate\write\gtoutfile{Authors: \ifx\theasciiauthors\relax
\theauthors\else\theasciiauthors\fi}
{\def\\{ }\immediate\write\gtoutfile{Title: \ifx\theasciititle\relax
\thetitle\else\theasciititle\fi}}
\immediate\write\gtoutfile{Subj-class: GT or SG, GR etc}
\immediate\write\gtoutfile{MSC-class: \theprimaryclass\ifx\thesecondaryclass\relax\else, \thesecondaryclass\fi}
\immediate\write\gtoutfile{Journal-ref: Algebr. Geom. Topol. \thevolumenumber\s
(\thevolumeyear) \startpage-\finishpage}
\immediate\write\gtoutfile{Comments: Published by Algebraic and
Geometric Topology at}
\immediate\write\gtoutfile{\s\s\s  http://www.maths.warwick.ac.uk/agt/AGTVol\thevolumenumber/agt-\thevolumenumber-\thepapernumber.abs.html}
\immediate\write\gtoutfile{\noexpand\\}
\immediate\write\gtoutfile{}
\ifx\theasciiabstract\relax
\immediate\write\gtoutfile{\theabstract}\else
\immediate\write\gtoutfile{\theasciiabstract}\fi
\immediate\write\gtoutfile{}
\immediate\write\gtoutfile{\noexpand\\}
\immediate\write\gtoutfile{}
\immediate\closeout\gtoutfile}}  

\def\maketitlepage{\makeagttitle\makeheadfile}
\let\makeshorttitle\maketitlepage
\let\maketitle\maketitlepage

\lognumber{51}
\volumenumber{4}
\volumeyear{2004}
\papernumber{51}
\pagenumbers{1155}{1175}
\received{6 May 2004} 
\revised{8 November 2004}
\accepted{15 November 2004}
\published{10 December 2004}

\usepackage[all]{xy}

\usepackage{amsmath}
\usepackage{amssymb}
\usepackage{stmaryrd}
\usepackage{epsfig}
\usepackage{bbm,color}
\usepackage{dbnsymb}

\newtheorem{thm}{Theorem}[section]
\newtheorem{lem}[thm]{Lemma}
\newtheorem{prop}[thm]{Proposition}
\newtheorem*{thm*}{Theorem}
\theoremstyle{definition}
\newtheorem{defn}[thm]{Definition}

\newtheorem*{quest}{Questions}

\newcommand{\sg}{\langle}
\newcommand{\sd}{\rangle}
\newcommand{\boA}{\mathcal{A}}
\newcommand{\boB}{\mathcal{B}}
\newcommand{\boD}{\mathcal{D}}

\newcommand{\boS}{\mathcal{S}}
\newcommand{\boC}{\mathcal{C}}

\newcommand{\Q}{\mathbb{Q}}
\newcommand{\Z}{\mathbb{Z}}
\newcommand{\rig}{\rightarrow}
\newcommand{\lef}{\leftarrow}
\newcommand{\alg}[3]{\vphantom{#2}_{#1} #2_{#3}}
\newcommand{\dd}{\text{d}}
\newcommand{\arete}{\!\!\frown}
\newcommand{\Zr}{Z^{\rm rat}}
\newcommand{\Zwr}{Z^{\fourwheel {\rm rat}}}
\newcommand{\disj}{\small{\amalg}}
\DeclareMathOperator{\lift}{Lift}
\DeclareMathOperator{\hair}{Hair}

\begin{document}
\title{A computation of the Kontsevich integral\\of torus knots}
\authors{Julien March\'e}
\coverauthors{Julien March\noexpand\'e}
\asciiauthors{Julien Marche}
\address{Institut de Math\'ematiques de Jussieu, \'Equipe ``Topologie et G\'eom\'etries Alg\'ebriques''\\
Case 7012, Universit\'e Paris VII, 75251 Paris CEDEX 05, France}
\asciiaddress{Institut de Mathematiques de Jussieu, Equipe 
`Topologie et Geometries Algebriques'\\
Case 7012, Universite Paris VII, 75251 Paris CEDEX 05, France}

\email{marche@math.jussieu.fr}
\begin{abstract}
We study the rational Kontsevich integral of torus knots.  We
construct explicitely a series of diagrams made of circles joined
together in a tree-like fashion and colored by some special rational
functions. We show that this series codes exactly the unwheeled
rational Kontsevich integral of torus knots, and that it behaves very
simply under branched coverings. Our proof is combinatorial. It uses
the results of Wheels and Wheeling and various spaces of diagrams.
\end{abstract}
\maketitle

\section{Introduction and notation}

In 1998, in \cite{bgrt}, a conjecture was formulated about a precise expression for the Kontsevich integral of the unknot (later it was proved, see for instance \cite{blt}). Until now, we do not know any complete formula for this powerful invariant for non-trivial
knots. Let $\boB$ be the usual space of uni-trivalent diagrams, which is the target space of the Kontsevich integral. We define a localization of $\boB$ called $\boB_s$. 
Unfortunately,
the localization map is not injective in high loop degrees.
 Then, we give a formula for the family of torus knots which takes its values in $\boB_s$ and which we state in a ``rational'' form.

The starting point of our proof is a well-known formula: it has been
used by Christine Lescop (see \cite{les}) and Dror Bar-Natan in
unpublished work.  From this expression, we construct a sequence of
series of diagrams which are all obtained by ``gluing wheels'' and
which converges to the unwheeled Kontsevich integral of torus knots.
(See Section 1.1 for the precise meaning of ``unwheeled'' here.)  From
this we present the result in a compact way by using ``gluing
graphs'', which we will define along the way.

\begin{figure}[htbp]
\qua \input{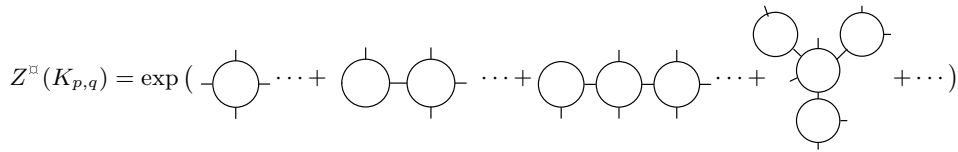}
\caption{Diagrams appearing in the unwheeled Kontsevich integral of 
torus knots}
\label{loop}
\end{figure}

 Then, in the third part, we compute a rational form of the preceding expression and show that only tree-like gluing diagrams appear as is suggested by figure \ref{loop}. More precisely, we show the following theorem:

\begin{thm*}
 Let $D$ be the operator on $\Q(t)$ defined by $Dg (t)=tg'(t)$ and let $h(t)=\frac{t+1}{t-1}$ (the operator $D$ acts as $\frac{\dd}{\dd x}$ on $g(\exp(x))$).

There is a  series of diagrams $Y^{rat}_{p,q}$ obtained by inserting circles in vertices of tree graphs, such that circles corresponding to vertices of valence $k$ are colored by $D^{k-1}h(t^p), D^{k-1} h(t^q)$ or $D^{k-1} h(t^{pq})$.

Applying the ``hair map'' (in other words, substituting $t$ with the exponential of a small leg attached to the circle),
we obtain a series in $\boB_s$ which is equal to the logarithm of the unwheeled Kontsevich integral of the torus knot of parameters $p,q$ plus a fixed series $\log \langle \Omega,\Omega \rangle$ (see 
\ref{intro} for precise definitions and theorem \ref{theoreme} for an explicit expression of the first terms).

\end{thm*}

As a consequence of this computation, we show that the operator $\lift_r$ which corresponds to cyclic branched coverings of $S^3$ along the knot simply acts on $Y_{p,q}^{rat}$ by multiplying a diagram $D$ by $r^{-\chi(D)}$ where $\chi$ is the Euler characteristic.

This article extends the result of our previous article \cite{moi} where we computed the unwheeled Kontsevich integral of torus knots up to degree 3. Here, we give a formula for all degrees using quite different methods.

Lev Rozansky also computed formulas for the loop expansion of torus knots in the weight system associated to $\rm{sl}_2$, see \cite{roz2}.
The computation of the 2-loop part of torus knots has been done independantly by Tomotada Ohtsuki in \cite{oht} who computed more generally a formula for 2-loop part of knots cabled by torus knots.

We would like to thank Stavros Garoufalidis, Marcos Marino, Tomotada Ohtsuki and Pierre Vogel for useful remarks. We also thank the referee and Gregor Masbaum for their remarks and their careful reading.

\subsection{Normalizations of the Kontsevich integral}\label{intro}

Let $K$ be a knot in $S^3$ and suppose that $K$ has a banded structure with self-linking 0. We will denote by $Z(K)$ the Kontsevich integral of $K$ in the algebra $\boA$ of trivalent diagrams lying on a circle.

Let $\boB$ be the algebra of uni-trivalent diagrams. It is well known that the Poincar\'e-Birkhoff-Witt map $\chi:\boB\to\boA$ is an isomorphism but not an algebra isomorphism. We will denote by $\sigma$ its inverse.

If $U$ is the trivial knot, we define $\Omega=\sigma Z(U)$. The series $\Omega$ is the exponential of a series of connected graphs whose first terms are $\frac{1}{48}\twowheel-\frac{1}{5760}\fourwheel+\cdots$.

The map $\Upsilon=\chi \circ \partial_{\Omega}:\boB\to \boA$ defined
for instance in \cite{th} is known to be an algebra isomorphism. The
quantity $Z^{\fourwheel}(K)=\Upsilon^{-1}Z(K)$ will be called
unwheeled Kontsevich integral. It behaves better than $\sigma Z(K)$
under connected sum and cyclic branched coverings.

For each knot $K$, the quantities $Z(K)$, $\sigma Z(K)$ and $Z^{\fourwheel}(K)$ are group-like, which means that they are exponentials of a series of connected diagrams. We will denote by respectively $z(K)$, $\sigma z(K)$ and $z^{\fourwheel}(K)$ the logarithms of these quantities.

\subsection{Loop degree and rationality}

If $D$ is a connected diagram of $\boB$, its first Betti number defines a degree called loop degree. The loop degree 1 part of $\sigma Z(K)$ or $Z^{\fourwheel}(K)$ is well-known: it only depends on the Alexander polynomial of $K$. For the higher degrees, very little is known. There are formulas for the 2-loop part of small knots in Rozansky's table (see \cite{roz}), and we can find in \cite{double} a formula for the 2-loop part of untwisted Whitehead doubles. In the sequel, we give a formula for the full Kontsevich integral of torus knots.

In order to make precise computations, we will need the following
formalism which is yet another version of ``diagrams with beads'' used
in \cite{roz, kr, rat} for instance, and generalized by P. Vogel in
\cite{vog}.

\begin{defn}
Let $\boC$ be the category whose objects are free abelian groups of finite rank and morphisms are linear isomorphisms. We call $\boC$-module a functor from $\boC$ to the category of $\Q$-vector spaces. In the sequel,
we will associate to any $\boC$-module $F$ a $\Q$-vector space $\boD(F)$ which will be called the space of diagrams decorated by $F$. This construction will be a functor from the category of $\boC$-modules to the category of $\Q$-vector spaces.
\end{defn}\label{espdediag}
Let $F$ be a $\boC$-module and $\Gamma$ be a finite trivalent graph with local orientations at vertices  (we allow $\Gamma$ to have connected components which are circles).

We define $\boD(F)$ as the quotient of $\bigoplus_{\Gamma} F(H^1(\Gamma,\Z))$ by the following relations:
\begin{itemize}
\item
If $\Gamma$ is isomorphic to $\Gamma'$ via a map $\phi$, then we identify $x\in F(\Gamma')$ and $F(\phi^*)(x)\in F(\Gamma)$ for any $x$ in $F(\Gamma')$.
\item
If $\Gamma$ and $\Gamma'$ just differ by the orientation of a vertex, then we identify $x$ in $F(H^1(\Gamma,\Z))$ with $-x$ in $F(H^1(\Gamma',\Z))$.
\item
If $\Gamma$ is a graph with one four-valent vertex, we note $\Gamma_I, \Gamma_H, \Gamma_X$ the three standard resolutions of this vertex. We have canonical identifications between $H^1(\Gamma,\Z), H^1(\Gamma_I,\Z),H^1(\Gamma_H,\Z)$ and $H^1(\Gamma_X,\Z)$. For every $x\in F(H^1(\Gamma,\Z))$, we add the relation $x_I=x_H-x_X$. 
\end{itemize}

Let us give the main example: we define $F(H)=\Q[[H]]=\prod_{n\ge 0}S^n(H\otimes \Q)$. The space $\boD(F)$ is obtained by coloring graphs with 1-cohomology classes which can be materialized by small legs attached to the edges. It is not hard to see that $\boD(F)\simeq\boB$. This isomorphism is a convenient way to express series of diagrams.

Let $f(x)$ be the power series defined by $\frac{1}{2}\log{\frac{\sinh{x/2}}{x/2}}$. The famous wheel formula (see \cite{th}) states that $\sigma z(U)=f(x)$ (we suppose that the variable $x$ is a generator of $H^1(\bigcirc,\Z)$). Further, it was shown in \cite{kr} that the loop degree 1 part of $\sigma z(K)$ is $f(x)+Wh_K(x)$ where $Wh_K(x)=-\frac{1}{2}\log \Delta(e^x)$ and $\Delta$ is the Alexander polynomial of $K$.
As we are interested in the higher loop degree part, we need to recall the rationality theorem which was proved in \cite{rat}. For that, we define a new space of diagrams with the help of a new functor.

Let $\Q[\exp(H)]$ be the polynomial algebra on elements on the form $\exp(\lambda), \lambda\in H$ with the relation $\exp(\lambda+\mu)=\exp(\lambda)\exp(\mu)$ for all $\lambda$ and $\mu$ in $H$.

Then, $F(H)=\Q[\exp(H)]_{loc}$ is the localization of the preceding space by expressions such as $P_1(\exp(h_1))\cdots P_k(\exp(h_k))$ with $h_1,\ldots,h_k \in H$, $P_1,\ldots,P_k\in \Q[h]$ and $P_1(1)\cdots P_k(1)\ne 0$. 

It defines a new space of diagrams $\boB^{\rm rat}=\boD(F)$ with a map $\hair:\boB^{\rm rat}\to\boB$ induced by the Taylor expansion: $\Q[\exp(H)]_{loc}\to\Q[[H]]$ which is a map of $\boC$-modules.

The rationality theorem tells us that there is an element $\Zr(K)\in \boB^{\rm rat}$ such that 
$$\sigma Z(K)=\exp(f(x)+Wh_K(x)) \hair \Zr(K).$$ 
Respectively, there is an element $\Zwr(K) \in \boB^{\rm rat}$ such that 
$$Z^{\fourwheel}(K)=
\frac{1}{\langle \Omega,\Omega\rangle}
\exp(f(x)+Wh_K(x)) \hair \Zwr(K).$$ 
A construction, developed in \cite{rat} gives $\Zr(K)$ as an invariant of $K$; this is not automatic because the $\hair$ map is not injective, 
although it is so in small degrees (see \cite{ber}). We refer to \cite{lift} for a construction of $\Zwr(K)$.

\section{Diagrammatic expressions of the integral}

\subsection{First diagrammatic expression}

Let $p$ and $q$ be two coprime integers such that $p>0$. We note $K_{p,q}$ the torus banded knot with parameters $p$ and $q$ and self-linking 0, and $L_{p,q}$ the torus banded knot with banding parallel to the torus on which it lies. This knot has self-linking $pq$, and its
Kontsevich integral is a bit easier to compute. 

The method of computation is inspired from \cite{les}: we first compute the Kontsevich integral of the following braid.

\begin{figure}[htbp]
\begin{center}
\input{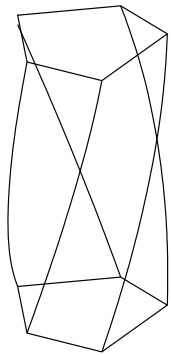}
\caption{Toric braid with parameters $(5,1)$}
\label{polygone}
\end{center}
\end{figure}

Let $p$ points be lying on the vertices of a regular $p$-gon. We note $\gamma$ the braid obtained by rotating the whole picture by an angle $2\pi\frac{q}{p}$ as is shown in figure \ref{polygone}.

Let us associate to any one dimensional manifold $\Gamma$ the space $\boA(\Gamma)$ of trivalent diagrams lying on $\Gamma$. This defines a contravariant functor with respect to continuous maps relative to boundaries. Let $\phi_p^*$ be the map induced by the projection onto
the first factor $\phi_p:[0,1]\times\{1,\ldots,p\}\to [0,1]$ and $\isolatedchord$ be the only degree 1 diagram in $\boA([0,1])$. 
Then a direct computation of monodromy of the K-Z connection shows that $\gamma$ has Kontsevich integral equal to $\phi_p^*( \exp_{\#}(\frac{q}{2p}\isolatedchord))$.

The banded knot $L_{p,q}$ is obtained by closing the previous braid: this translates diagrammatically to the following: let $\psi_p$ be the map from $S^1$ to itself defined by $\psi_p(z)=z^p$. Then, $Z(L_{p,q})=\psi_p^*(\nu\# \exp_{\#}(\frac{q}{2p}\isolatedchord))$, where $\nu=Z(U)$.

By lemma 4.10 of \cite{th}, the map $\psi_p^*$ viewed in $\boB$ has the following form:
if $D\in \boB$ has $k$ legs (i.e. univalent vertices) then $\sigma \psi_p^* \chi D  = p^k D$. We will note more simply $D_p$ the result of this operation which looks like a change of variable.

Then, to compute $Z(K_{p,q})$ from $Z(L_{p,q})$, we only need to change the framing, that is $Z(K_{p,q})=\exp_{\#}(-\frac{pq}{2}\isolatedchord)\#Z(L_{p,q})$. We will transform this product into
the usual one by applying the unwheeling map $\Upsilon^{-1}$. As a result, we will have a formula for $Z^{\fourwheel}(K_{p,q})$.

We now sum up the steps of the computation:

\begin{enumerate}
\item Computation of $\sigma(\nu\# \exp_{\#}(\frac{q}{2p}\isolatedchord))$
\item Change of variables $x\mapsto px$
\item Unwheeling
\end{enumerate}

We recall that $\Upsilon=\chi\circ\partial_{\Omega}$ is an algebra isomorphism and that $\Upsilon^{-1}\nu = \frac{\Omega}{\langle \Omega,\Omega\rangle}$ and $\Upsilon^{-1}\isolatedchord= \frown-\ThetaGraph/24$. Here $\langle A,B\rangle$ is the sum over all ways of gluing all legs of $A$ to all legs of $B$.

Then, to achieve the first step, we have to compute
\begin{equation}\label{expr}
\sigma(\nu\# \exp_{\#}(\frac{q}{2p}\isolatedchord))=
\frac{\partial_{\Omega} (\Omega\exp(\frac{q}{2p}\arete))}{\langle\Omega,\Omega\rangle\exp(\frac{q}{48p}\ThetaGraph)}.
\end{equation}
We explain below what will be a diagram colored by some parameters. In this setting, we will recall how usual operations on diagrams are obtained. Finally we will describe a new operation.

\begin{defn}\label{definition}$\phantom{99}$
\begin{itemize}
\item
Let $P$ be a set of parameters. We denote
$\boB(P)$ the space of couples $(D,f)$ where $D$ is a uni-trivalent graph and $f$ is a map from legs of $D$ to $P$. We will say that legs of $D$ are labeled or colored by elements of $P$. If $D\in\boB$ and $x\in P$, we will write $D_x$ the diagram $D$  whose legs are colored by $x$. 

\item If $D\in \boB(x)$, we define the diagram $D_{x+y}\in\boB(x,y)$ by replacing legs of $D$ by the same legs colored by $x$ or $y$ in all possible ways.
\item If $D,E\in\boB(P)$ and $x\in P$ we define $\sg D,E \sd_x$ as the sum over all gluings of all $x$-legs of $D$ on all $x$-legs of $E$.
We also define $\partial_D E$ as the sum over all gluings of all $x$-legs of $D$ on some $x$-legs of $E$.
This operator satisfies $\partial_{D_x} E_x = \sg D_y,E_{x+y}\sd_y$ and  for $F\in\boB(P)$, $\sg D_x,E_x F_x\sd_x=\sg\partial_{E_x}D_x,F_x\sd_x$.

\item
Let $A_x$ and $B_x$ be two series of diagrams in $\boB(P)$ where $P$ contains at least three colors: $\{x,y,z\}$. We define $A\cdot B=\langle A_{y+x},B_{x+z}\rangle_x$. The $y$-legs of $A\cdot B$ will be called left legs and $z$-legs of $A\cdot B$ will be called right legs for obvious reasons. For rational numbers $r$ and $r'$, the diagram $\alg{r}{A\cdot B}{r'}$ will be the result of multiplying all left legs by $r$ and all right legs by $r'$.

Moreover, if $a_x$ and $b_x$ are two series of connected diagrams, then the series $\exp(a), \exp(b)$ and $\exp(a)\cdot\exp(b)$ are group-like, hence, we define  $\alg{y}{a\times b}{z}=\log(\alg{y}{\exp(a)\cdot\exp(b)}{z})$.
\end{itemize}
\end{defn}

Using the notation of definition \ref{definition}, we compute
$$\partial_{\Omega}(\Omega\exp(\frac{q}{2p}\arete))=\langle \Omega_x,\Omega_{x+y}\exp(\frac{q}{2p}\alg{x+y}{\arete}{x+y})\rangle_x.$$
We use the fact that $\alg{x+y}{\arete}{x+y}=\alg{x}{\arete}{x}+2\alg{x}{\arete}{y}+\alg{y}{\arete}{y}$ and properties recalled in the preceding definition 
to obtain that the series $\partial_{\Omega}(\Omega\exp(\frac{q}{2p}\arete))$
is equal to
$\langle \partial_{\exp(\frac{q}{2p}\alg{x}{\arete}{x})}\Omega_x,\Omega_{x+y}\exp(\frac{q}{p}\alg{x}{\arete}{y})\rangle_x \exp(\frac{q}{2p}\alg{y}{\arete}{y})$.

Now we use the fundamental formula $\partial_{D_x}\Omega_x = \sg D_x,\Omega_x\sd_x \Omega_x$ which is true if $D$ contains only $x$-legs.
The preceding expression reduces to $$\langle \exp(\frac{q}{2p}\arete),\Omega\rangle \left(\alg{\frac{q}{p}}{\Omega\cdot\Omega}{}\right)\exp(\frac{q}{2p}\arete).$$
Thanks to the identity $\langle \exp(\frac{q}{2p}\arete),\Omega\rangle=\exp(\langle \frac{q}{2p}\arete,\Omega\rangle)=\exp(\frac{q}{48p}\ThetaGraph)$, we 
may cancel this factor in the expression \eqref{expr} above.

Multiplying by $p$ to the power the number of legs, we obtain a formula for $Z(L_{p,q})$:
$$Z(L_{p,q})= \frac{\alg{q}{\Omega\cdot\Omega}{p}\exp(\frac{qp}{2}\arete)}{\langle\Omega,\Omega\rangle}.$$
We obtain a formula for $Z^{\fourwheel}(K_{p,q})$ by unwheeling $Z(L_{p,q})$ and multiplying by $\Upsilon^{-1} \exp_{\#}(-\frac{pq}{2}\isolatedchord)=\exp(-\frac{pq}{2}\arete+\frac{pq}{48}\ThetaGraph)$:
\begin{equation}\label{fff}
Z^{\fourwheel}(K_{p,q})=
\partial_{\Omega}^{-1}(
\alg{q}{\Omega\cdot\Omega}{p}\exp(\frac{qp}{2}\arete)
)
\frac{\exp(-\frac{pq}{2}\arete+\frac{pq}{48}\ThetaGraph)}{\langle\Omega,\Omega\rangle}.
\end{equation}

\subsection{A sequence converging to the integral of torus knots}

In this part, we will start from formula \eqref{fff} to deduce new formulas by an induction process.

We define  $\boB^c= \boB(\{\text{active,inert}\})$. There is a forgetful
map $\boB^c\to\boB$.

Let us define three operators analogous to those of definition \ref{definition} but for diagrams in $\boB^c$ in the following way:

\begin{defn}\label{definition2}$\phantom{99}$
\begin{itemize}
\item
For $A\in \boB^c$ and $r\in\Q$, we set $A_r$ to be the diagram $A$ whose active legs are multiplied by $r$, viewed as an element of $\boB$.
\item
If $A\in\boB^c$ and $B\in\boB$, the diagram $\partial_A B \in\boB$ is the sum over gluings of all active legs of $A$ on some legs of $B$. This is an element of $\boB$.
\item
For $A\in\boB$ and $B\in \boB^c$ two diagrams, we set $A\cdot B\in\boB^c$ to be the sum over all gluings of some active legs of $B$ with some legs of $A$. The active legs of $A\cdot B$ are the remaining active legs of $B$. We also define $a\times b=\log(\exp(a)\times\exp(b))$. See figure \ref{active} for both definitions.
\end{itemize}
\begin{figure}[htbp]
\begin{center}
\input{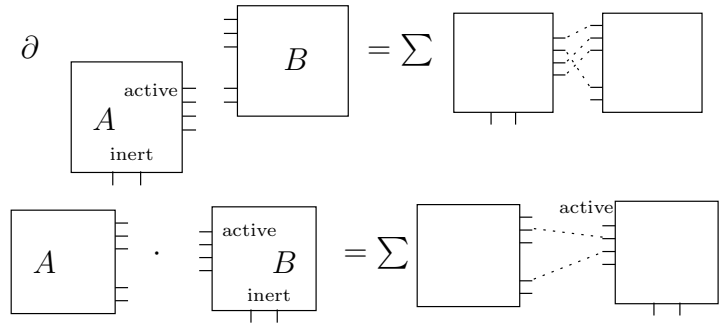}
\caption{The operators $\partial$ and $\cdot$}
\label{active}
\end{center}
\end{figure}
\end{defn}

\begin{lem}\label{lm}
For all $A\in \boB^c$ and $B,C$ in $\boB$ we have the identity $\partial_A(BC)=\partial_{C\cdot A}B$.
\end{lem}

\begin{figure}[htbp]
\begin{center}
\input{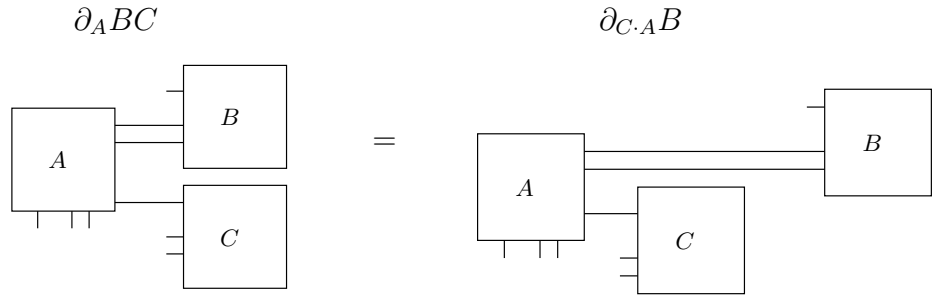}
\caption{Graphical proof of lemma \ref{lm}}
\label{lmfig}
\end{center}
\end{figure}

\begin{proof}
In the formula $\partial_A(BC)$, the active legs of $A$ are glued to legs of $B$ and $C$. We can split these legs in two parts and glue the first part to $C$ and the second part to $B$. Before gluing the second part to $B$, we obtain exactly the diagram $C\cdot A$. The remaining active legs of $A$ should be glued to $B$: we write the result of this operation $\partial_{C\cdot A}B$.
This proof is summed up in figure \ref{lmfig}.\end{proof}

We now define a sequence of series $\omega^n$ in $\boB^c$.
Let us put $\omega^{-1}=\alg{\frac{1}{p}}{\omega\times\omega}{\frac{1}{q}}$ and $\omega^0=\omega$. We declare that all their legs are active.

We suppose all legs of $\omega^{-1}$ and $\omega^0$ are active, and for all $n\ge 0$, we set $\omega^{n+1}=(\omega^{n-1}_{pq}-\omega^{n}_{pq})\times\omega^n -\omega^{n-1}_{pq}+\omega^{n}_{pq}$. In this formula, the added term $\omega^{n-1}_{pq}-\omega^{n}_{pq}$ is considered inert.

Let us prove the following proposition:

\begin{prop}
For all $n\ge 0$, we have
\begin{enumerate}
\item[\rm(i)]
The series $\exp(\frac{pq}{2}\arete)$ acts by derivation on active legs of $\exp(\omega^n)$ by multiplication by $\exp(\frac{pq}{48}\ThetaGraph)$. In formulas, we mean that $\partial^{\text{active}}_{\exp(\frac{pq}{2}\arete)}\exp(\omega^n)=\exp(\frac{pq}{48}\ThetaGraph)\exp(\omega^n)$. Here we have written $\partial^{\text{active}}$ in order to mean that only active legs are glued in this process, contrary to the definition \ref{definition2}.

\item[\rm(ii)]$\phantom{99}$\vspace{-15pt}
$$Z^{\fourwheel}(K_{p,q})=\left(\partial^{-1}_{\exp(\omega^n)}\exp(\frac{pq}{2}\arete+\omega^{n-1}_{pq})\right)\frac{\exp(-\frac{pq}{2}\arete\!+\!\omega^{-1}_{pq}\!-\!\omega^{n-1}_{pq}\!+\!\frac{pq}{48}\ThetaGraph)}{\langle\Omega,\Omega\rangle}$$
\item[\rm(iii)]$\phantom{99}$\vspace{-15pt}
\begin{multline}
Z^{\fourwheel}(K_{p,q})=\partial^{-1}_{\exp(\omega^n)}\left[\partial_{\exp(\omega^n)}\exp(\frac{pq}{2}\arete)\exp(\omega^{n-1}_{pq}-\omega^n_{pq})\right]\cdot\\
\frac{\exp(-\frac{pq}{2}\arete+\omega^{-1}_{pq}-\omega^{n-1}_{pq})}{\langle\Omega,\Omega\rangle}.
\end{multline}

\item[\rm(iv)]

Moreover the sequence $\omega^{-1}_{pq}-\omega^n_{pq}$ converges to $\log\left(Z^{\fourwheel}(K_{p,q})\langle\Omega,\Omega\rangle\right)$ with respect to loop degree.
\end{enumerate}
\end{prop}

\begin{proof}
(i)\qua
Let us prove the first formula by induction.

The assertion for $n=0$ is just $\partial_{\exp(\frac{pq}{2}\arete)}\Omega=\exp(\frac{pq}{48}\ThetaGraph) \Omega$. This is a consequence of the fundamental relation $\partial_D\Omega = \sg D,\Omega\sd \Omega$.

Suppose that for some $n\ge 0$ we have the identity $$\partial^{\text{active}}_{\exp(\frac{pq}{2}\arete)}\exp(\omega^n)=\exp(\frac{pq}{48}\ThetaGraph)\exp(\omega^n).$$ 
Then by definition, we have the identity $$\exp(\omega^{n+1})=\frac{\exp(\omega^{n-1}_{pq}-\omega^n_{pq})\cdot \exp(\omega^n)}{\exp(\omega^{n-1}_{pq}-\omega^n_{pq})}.$$
But all active legs of the series $\exp(\omega^{n+1})$ come from the series $\exp(\omega^n)$. This proves the identity $$\partial^{\text{active}}_{\exp(\frac{pq}{2}\arete)}\exp(\omega^{n+1})=\frac{\exp(\omega^{n-1}_{pq}-\omega^n_{pq})\cdot \partial^{\text{active}}_{\exp(\frac{pq}{2}\arete)}\exp(\omega^n)}{\exp(\omega^{n-1}_{pq}-\omega^n_{pq})}$$ because all ``derivation-like'' operations commute (see for instance lemma \ref{lm}).

 Using the induction assumption, we finally prove the following formula:
$$\partial^{\text{active}}_{\exp(\frac{pq}{2}\arete)}\exp(\omega^{n+1})=\exp(\frac{pq}{48}\ThetaGraph) \exp(\omega^{n+1}).$$

(ii)\qua
The second formula for $n=0$ is nothing but a version of formula \eqref{fff}.
Suppose the second formula is true for some $n\ge 0$. To prove the second formula for $n+1$, we will prove that the expression at rank $n+1$ is equal to the expression at rank $n$.
This is equivalent to the following identity:
\begin{multline*}
\partial^{-1}_{\exp(\omega^{n+1})}
\exp(\frac{pq}{2}\arete+\omega^{n}_{pq})\exp(-\omega^n_{pq})
=\\
\partial^{-1}_{\exp(\omega^n)}
(\exp(\frac{pq}{2}\arete+\omega^{n-1}_{pq}))\exp(-\omega^{n-1}_{pq}).
\end{multline*}
By applying the operator $\partial_{\exp(\omega^{n+1})}$, this identity is equivalent to the following equation:
$$
\exp(\frac{pq}{2}\arete+\omega^{n}_{pq})=\partial_{\exp(\omega^{n+1})}[
\partial^{-1}_{\exp(\omega^n)}
(\exp(\frac{pq}{2}\arete+\omega^{n-1}_{pq}))\exp(\omega^{n}_{pq}-\omega^{n-1}_{pq})
]$$
Using lemma \ref{lm}, we have the following equalities:
\begin{eqnarray*}
\partial_{\exp(\omega^{n+1})} A 
&=& \partial_{\frac{\exp(\omega^{n-1}_{pq}-\omega^n_{pq})\cdot \exp(\omega^n)}{\exp(\omega^{n-1}_{pq}-\omega^n_{pq})}}A \\
&=& \left(\partial_{\exp(\omega^{n-1}_{pq}-\omega^n_{pq})\cdot \exp(\omega^n)} A \right) \exp(\omega^n_{pq}-\omega^{n-1}_{pq})\\
&=& \partial_{\exp(\omega^n)}\left(A\exp(\omega^{n-1}_{pq}-\omega^n_{pq})\right) \exp(\omega^n_{pq}-\omega^{n-1}_{pq})
\end{eqnarray*}
Replacing $A$ by its value, we obtain the desired identity.

(iii)\qua
We obtain the third formula from the second one by using the following formula:
 \begin{eqnarray*}
\partial_{\exp(\omega^n)}\exp(\frac{pq}{2}\arete)&=&\langle \exp(\omega^n_a),\exp(\frac{pq}{2}\alg{a+b}{\arete}{a+b})\rangle_a\\
&=& \left(\partial^{\text{active}}_{\exp(\frac{pq}{2}\arete)} \exp(\omega^n)\right)_{pq}\exp(\frac{pq}{2}\arete) \\
&=& \exp(\omega^n_{pq})\exp(\frac{pq}{48}\ThetaGraph+\frac{pq}{2}\arete)
\end{eqnarray*}
The last step uses the first part of the proposition.

(iv)\qua
Let us prove the last property: we will use mainly the third formula of the proposition.

We want to show that all diagrams in the series 
\begin{equation}\label{rappel}
\partial^{-1}_{\exp(\omega^n)}\left[\partial_{\exp(\omega^n)}\exp(\frac{pq}{2}\arete)\exp(\omega^{n-1}_{pq}-\omega^n_{pq})\right]\exp(-\frac{pq}{2}\arete).
\end{equation}
have loop degree greater or equal to $n$. This will prove that the expression $\omega^{-1}_{pq}-\omega^{n-1}_{pq}$ and the expression $\log(Z^{\fourwheel}(K_{p,q})\langle\Omega,\Omega\rangle)$ coincide up to loop degree $n$, and then that the series $\omega^{-1}_{pq}-\omega^{n-1}_{pq}$ converges to  $\log(Z^{\fourwheel}(K_{p,q})\langle\Omega,\Omega\rangle)$ with respect to loop degree when $n$ goes to infinity.

Let us prove by induction that for all $n$, the series $\omega^{n-1}_{pq}-\omega^n_{pq}$ contains diagrams with loop degree greater or equal to $n$. 
This is obvious for $n=0$. Suppose it is true for some $n$, then by definition of $\omega^{n+1}$, we have 
$$\omega^{n+1}_{pq} -\omega^n_{pq}= \left((\omega^{n-1}_{pq}-\omega^n_{pq})\times \omega^n\right)_{pq}-\omega^{n-1}_{pq}.$$
This proves that the series $\omega^{n+1}_{pq} -\omega^n_{pq}$ is obtained by non trivial gluings of $\omega^{n-1}_{pq}-\omega^n_{pq}$ and $\omega^n$. As by assumption the series $\omega^{n-1}_{pq}-\omega^n_{pq}$ has loop degree greater or equal to $n$, the series $\omega^{n+1}_{pq} -\omega^n_{pq}$ must have loop degree greater than $n$.

Now, the operator $\partial_{\exp(\omega^n)}$ always increases the loop degree, and so does its inverse. 
This proves that diagrams of $\exp(\omega^{n-1}_{pq}-\omega^n_{pq})$ cannot create  diagrams of loop degree less than $n$ in the expression \eqref{rappel}.
 Hence, up to degree $n$, we can ignore the series $\exp(\omega^{n-1}_{pq}-\omega^n_{pq})$, and the preceding expression reduces to 
$$\partial^{-1}_{\exp(\omega^n)}\left(\partial_{\exp(\omega^n)}\exp(\frac{pq}{2}\arete)\right)\exp(-\frac{pq}{2}\arete)=1.$$ 
This proves the last part of the proposition.\end{proof}

\subsection{An expression with gluing graphs}

We remark that all diagrams appearing in $\omega^n$ are made by ``gluing wheels'', i.e. they are constructed from wheel series as $\omega,\omega_p,\omega_q$ by applying some gluing operators. We make this statement precise
by introducing gluing graphs and we give a presentation of $\log(Z^{\fourwheel}(K_{p,q})\langle\Omega,\Omega\rangle)$ with gluing graphs. These graphs will be useful for finding a rational expression of the preceding expression. As a non-trivial
application, we will show that only tree-like graphs appear in this expression.

\begin{defn}
Let $P$ be a set of parameters.
We denote by $\boS(P)$ the $\Q$-vector space generated by unoriented finite graphs $(V,E,h)$ with a map $V\to P$ where $V$ is the set of vertices and $E$ the set of edges and a map $h:P\to\Q[[x]]$. Elements of this space will be called gluing graphs colored by $P$. The number of edges is a degree on $\boS(P)$. We complete $\boS(P)$ with respect to this degree. This space has an obvious Hopf algebra structure whose primitive elements are formed of connected graphs.

We define a substitution map $s:\boS(P)\to \boB$ in the following way:
if $(X,h)$ is a diagram in $\boS(P)$, we define $s(X,h)$ 
 by gluing for all $a\in P$ and vertices $v$ decorated by $a$ the edges adjacent to $v$ to wheels generated by $h(a)$ in all possible ways. Sometimes, we call $a$-colored the remaining legs generated by $h(a)$. It is clear that $s$ is a map of Hopf algebras.
\end{defn}

Let us give some examples:
if $P=\{a\}$ $X=\bullet_a$ and $h(a)=f(x)$, then $s(X,h)=\omega$. In the same way, $s(\exp(\bullet_a),h)=\Omega$ and 
$$s({}_a\bullet\!\!-\!\!\bullet_b\!\!\!\!\!-\!\!\!-\!\!\!\bullet_c,h)=16 -\!\!\bigcirc\!\!-\!\!\bigcirc\!\!-\!\!\bigcirc\!\!\equiv$$
 where $h(a)=x^2, h(b)=x^2$ and $h(c)=x^4$.

Let us define two operations on $\boS(P)$:
\begin{defn}\label{operations}$\phantom{99}$
\begin{itemize}
\item
Given two connected gluing graphs $X$ and $Y$ with $A$ some set of parameters of $X$ and $B$ some set of parameters of $Y$, we define $X\alg{A}{\times}{B}Y=\log(\exp(X)\cdot\exp(Y))$ where $\exp(X)\cdot\exp(Y)$ is obtained by adding in all ways a finite number of edges from $A$-colored vertices of $\exp(X)$ to $B$-colored vertices of $\exp(Y)$.
\item
Let $(X,h)$ be a diagram in $\boS(P\disj \{a\})$, $h(x)\in \Q[[x]]$, and $r\in\Q\setminus\{0\}$.
We note $a^{r}$ the operator which divides
$X$ by $r^N$ where $N$ is the sum of valences of $a$-colored vertices of $X$.
\end{itemize}
\end{defn}

Let us look at
an example: take two formal parameters $a$ and $b$, then 
$$\exp(\bullet_a)\cdot\exp(\bullet_b)=\exp(\bullet_a+\bullet_b+{}_a\bullet\!\!-\!\!\bullet_b+{}_a\bullet\!\!\!=\!\!\!\bullet_b
+\frac{1}{2}{}_a\bullet\!\!-\!\!\bullet_b\!\!\!\!\!-\!\!\!-\!\!\!\bullet_a
+\frac{1}{2}{}_b\bullet\!\!-\!\!\bullet_a\!\!\!\!\!-\!\!\!-\!\!\!\bullet_b
+{}_a\bullet\!\!\!\equiv\!\!\!\bullet_b+ \cdots).$$
The operations of definition \ref{operations} are a version for gluing graphs of usual operations on diagrams as is shown in the next proposition.

\begin{prop}\label{compat}
Let $X$ and $Y$ be two gluing graphs in $\boS(P)$. Let $A$ and $B$ be two disjoint subsets of $P$. Then, we have $s(X\alg{A}{\times}{B}Y)=s(X)\alg{A}{\times}{B}s(Y)$ (we omit the decorations of the vertices are they remain the same).

Let $X$ be a graph in $\boS(P\disj\{a\})$. We have $s(X,h)_{ra}=s(a^r X, h(rx))$ (in this formula, we only change the decoration of the $a$-vertex).
\end{prop}
\begin{proof}
The proof of this proposition comes from a direct combinatorial description of the gluing operations that are involved, hence we omit it.
\end{proof}

The main result of this section is the following:

\begin{prop}\label{Xpq}
There is an explicit series of gluing graphs $X_{p,q}$ whose substitution is $\log(Z^{\fourwheel}(K_{p,q})\langle \Omega,\Omega \rangle)$.
\end{prop}
\begin{proof}
We will show by induction that for all $n\ge 0$, $\omega^n=s(X^n)$ for some $X^n\in\boS(P)$ where $P=\{*\}\cup\{a,b,c\}$.
The first parameter is active and corresponds to $f(x)$, the three last parameters are inert and correspond to $f(px),f(qx)$ and $f(pqx)$.

Firstly, the diagram $\omega^{-1}_{pq}$ has only inert legs and is obtained from the substitution of the gluing graph: $X^{-1}_{pq}=a^p b^q \bullet_a \times \bullet_b \in \boS(a,b)$ where $a$ and $b$ are respectively associated to $f(px)$ and $f(qx)$. We deduce directly this construction from the formula $\omega^{-1}_{pq}=\alg{p}{\omega\times\omega}{q}$ and proposition \ref{compat}.

We start the recursion at $n=0$ by setting $X^0=\bullet_*$. We have $s(X^0)=\omega^0$.

Take $n\ge 0$ and suppose we have constructed $X^k$ for all $k\le n$. 
Then, we set $X_{pq}^{k}= *^{pq} X^k |_{*\to c}$ for all $0\le k \le n$ in such a way that we have $\omega^k_{pq}=s(X_{pq}^k)$ for all $-1\le k \le n$ where $c$ is associated to $f(pqx)$.

Then we set $X^{n+1}=(X^{n-1}_{pq}-X^n_{pq})\alg{a,b,c}{\times}{*}X^n-(X^{n-1}_{pq}-X^{n}_{pq})$.
This definition satisfies $\omega^{n+1}=s(X^{n+1})\in \boB^c$.
It gives a recursive way for finding all diagrams $X^n$ and then we set $X_{p,q}=\lim_n (X^{-1}_{pq}-X^n_{pq})$.
\end{proof}

Let us give all terms with less than two edges as an example. We use some notations to give a more compact form. Edges may be oriented or not, and vertices are expressed as a sum of integers in a box.
To obtain the result, we color a vertex $\fbox{m+n}$ by $f(mx)+f(nx)$. Then, we obtain a sum of diagrams colored by integers. We divide each diagram of vertices $x_1,\ldots,x_n$ by the product of colors of each $x_i$ to the power the number of adjacent edges which are not coming to this vertex.

For instance the diagram $\fbox{pq}\!\!\rig\!\!\fbox{p}\!\!\lef\!\!\fbox{q}$ is a graphical expression for $\frac{1}{pq^2}{}_c\bullet\!\!-\!\!\bullet_a\!\!\!\!-\!\!-\!\!\!\bullet_b$ and the dots mean terms with more than two edges.
\begin{align*}
\omega^{-1}_{pq}&=\fbox{p}+\fbox{q}+\fbox{p}\!\!-\!\!\fbox{q}+\fbox{p}\!\!=\!\!\fbox{q}+\frac{1}{2}\fbox{p}\!\!-\!\!\fbox{q}\!\!-\!\!\fbox{p}
+\frac{1}{2}\fbox{q}\!\!-\!\!\fbox{p}\!\!-\!\!\fbox{q}+\cdots\\
\omega^0&=\fbox{1}\\
\omega^1&=
\fbox{1}+\fbox{p+q-pq}\!\!\lef\!\!\fbox{1}+
\fbox{p+q-pq}\!\!\leftleftarrows\!\!\fbox{1}+\frac{1}{2}\fbox{1}\!\!\rig\!\!\fbox{p+q-pq}\!\!\lef\!\!\fbox{1}\\
&+\frac{1}{2}\fbox{p+q-pq}\!\!\lef\!\!\fbox{1}\!\!\rig\!\!\fbox{p+q-pq}
+\fbox{1}\!\!\rig\!\!\fbox{p}\!\!-\!\!\fbox{q}+\fbox{1}\!\!\rig\!\!\fbox{q}\!\!-\!\!\fbox{p}+\cdots\\
\omega^2&=\omega^1 - \fbox{p+q-pq}\!\!\lef\!\!\fbox{pq}\!\!\lef\!\!\fbox{1} - \fbox{1}\!\!\rig\!\!\fbox{p+q-pq}\!\!\lef\!\!\fbox{pq}+\cdots\\
X_{p,q}&=\lim_{n}(\omega^{-1}_{pq}-\omega^n_{pq})=\fbox{p}+\fbox{q}-\fbox{pq}+\fbox{p}\!\!-\!\!\fbox{q} -\fbox{p+q-pq}\!\!\lef\!\!\fbox{pq}+
\fbox{p}\!\!=\!\!\fbox{q}\\
&-\fbox{p+q-pq}\!\!\leftleftarrows\!\!\fbox{pq}
+\frac{1}{2}\fbox{p}\!\!-\!\!\fbox{q}\!\!-\!\!\fbox{p}
+\frac{1}{2}\fbox{q}\!\!-\!\!\fbox{p}\!\!-\!\!\fbox{q}\\
&+\frac{1}{2}\fbox{pq}\!\!\rig\!\!\fbox{p+q-pq}\!\!\lef\!\!\fbox{pq}
-\frac{1}{2}\fbox{p+q-pq}\!\!\lef\!\!\fbox{pq}\!\!\rig\!\!\fbox{p+q-pq}\\
&-\fbox{pq}\!\!\rig\!\!\fbox{p}\!\!-\!\!\fbox{q}-\fbox{pq}\!\!\rig\!\!\fbox{q}\!\!-\!\!\fbox{p}
+\fbox{p+q-pq}\!\!\lef\!\!\fbox{pq}\!\!\lef\!\!\fbox{pq} +\cdots
\end{align*}

\section{Rationality}

The expression we have obtained until now is not rational. When trying to find a rational expression, we will show that expressions coming from non-tree gluing graphs vanish. For example, the graph $\fbox{p}\!\!=\!\!\fbox{q}$ above will not appear in the final
expression of the unwheeled Kontsevich integral of torus knots.

\subsection{Diagrams with singular colorings}\label{diag}

We need to work with diagrams which are colored by some singular algebraic expressions.
For that we define two spaces of diagrams $\boB_s$ and $\boB_s^{rat}$ which fit in the following commutative diagram:
$$\xymatrix{
\boB^{\rm rat}\ar[d]^{\hair}\ar[r] & \boB^{\rm rat}_s \ar[d]^{\hair} \\
\boB\ar[r] & \boB_s
}$$
Here we define $\boB_s=\boD(F)$ where $F(H)=(H\setminus\{0\})^{-1} \Q[[H]]$
, the localization of $\Q[[H]]$ by elements of $H\setminus\{0\}$. The space $\boB^{rat}_s$ is constructed from the functor $F(H)=S^{-1}\Q[\exp(H)]$, where $S$ is made of non zero expressions of the form $P_1(\exp(h_1))\cdots P_k(\exp(h_k))$ for $P_1,\ldots,P_k\in\Q[h]$ and $h_1,\ldots,h_k\in H$.

The $\hair$ map is again defined by the map induced by the 
Taylor expansion.

These constructions motivate the definition \ref{espdediag}. In the article \cite{moi}, we used a different construction using coloring of edges by elements of $\Q[[h]][h^{-1}]$. This construction was not convenient because the map from non-singular diagrams to singular ones fails to be injective even for small degrees as is shown in the following identity:
$$\xymatrix{
-\!\!\bigcirc\!\!-\!\!\bigcirc\!\!- = 
-\!\!\bigcirc\!\!\raisebox{2pt}{$\overset{\frac{1}{h} h}{\perp\!\!\perp}$}\!\!\bigcirc\!\!- = 0.
}$$
To avoid this, we have changed the construction of singular diagrams in some way which avoids
inverting
null-homologous legs.

\subsection{A formula for the substitution map}

Let $X$ be a connected gluing graph, with vertices $x_1,\ldots,x_N$ and power series $f_1,\ldots,f_N$ corresponding to these vertices. The aim of this part is to give an explicit formula for $s(X)\in\boB_s$.
We describe it in the following proposition:

\begin{prop}\label{subst}
Given a gluing graph $X$ with power series $f_1,\ldots,f_N$, its substitution is a finite combination of diagrams $\Gamma$ obtained by gluing at each vertex the incoming edges to a circle, and decorating by expressions of the form 
$$\prod\limits_{i=1}^{N}\left(\frac{f_i'(y_i)}{{y_i}^{p_i}}\right)^{(k_i-p_i)}/D_i$$
 where the $y_i$s are non zero cohomology classes, $k_i$ is the valency of the $i$-th vertex, $p_i$ is an integer satisfying $0\le p_i\le k_i$ and $D_i$ is a product of $p_i$ linear terms. Moreover, we have the inequality $\sum_i p_i>0$ unless $X$ is a tree. 
\end{prop}

\begin{proof}
To compute $s(X)$, we may glue at each vertex $x$ edges incoming to this vertex to wheels generated by $f(x)$. Hence, we can choose a cycling ordering of edges around each vertex, compute the result of the gluings which repect this order, and sum over such orderings.

Fix an ordering $e^x_1,\ldots,e^x_{k_x}$ of edges around each vertex $x$.
Let $\Gamma$ be the trivalent graph obtained by replacing $x$ by a circle attached to the edges in the prescribed order, and $H=H^1(\Gamma,\Z)$.
Fix the order $n_x$ of the wheel of $f(x)$ we will glue at $x$.
Write $x_1,\ldots,x_{k_x}$ the elements of $H$ coming from the edges $e_1$ to $e_2$, ..., $e_{k_x}$ to $e_1$ (see figure \ref{subs}).

\begin{figure}[htbp]
\begin{center}
\input{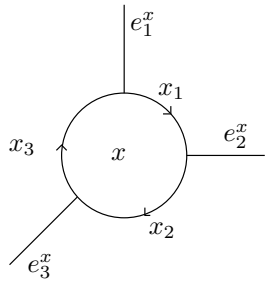}
\caption{Diagram obtained by substitution}
\label{subs}
\end{center}
\end{figure}

First, we glue the edge $e_1$ to any leg, so we have a factor $n_{x}$.
Then all gluings are listed in the following formula:
\begin{equation}\label{ratexp}
\sum\limits_{\text{orderings},x}n_x\sum\limits_{i_1+\cdots+i_{k_x}=n_x-k_x}x_1^{i_1}\cdots x_{k_x}^{i_{k_x}}.
\end{equation}
We explain an algorithmic reduction for this formula.
Suppose there are two indices $l$ and $m$ such that $x_l\ne x_m \in H^1(\Gamma,\Z)$, write
$i_{lm}=i_l+i_m$. Then, $\sum\limits_{i_1+\cdots+i_{k_x}=n_x-k_x}x_1^{i_1}\cdots x_{k_x}^{i_{k_x}}$ may be replaced by the sum 
$$\sum\limits_{i_1+\cdots+i_{k_x}+i_{lm}=n_x-k_x}x_1^{i_1}\cdots\widehat{x_l^{i_l}}\cdots \widehat{x_m^{i_m}} \cdots x_{k_x}^{i_{k_x}} \frac{x_l^{i_{lm}}-x_m^{i_{lm}}}{x_l-x_m} \in (H\setminus\{0\})^{-1}\Q[[H]]$$
This 
is a difference of 
sums of the same form but 
where each summand is now a product of $k_x-1$ monomials,
and with a degree one denominator.
Observe that the sum of the exponents of the monomials in each summand is 
unchanged and equal to $n_x-k_x$ as before.
We can perform such a reduction again for each 
one of the sums thus obtained 
a different number of times. This algorithm stops if all cohomology classes involved in 
the resulting expression
are the same.

Suppose that 
an expression has been obtained form \eqref{ratexp} by reducing it $p_x$ 
times. Then this 
expression
is a sum of identical 
terms each of which is
equal to a cohomology class $y$ to the power 
$n_x-k_x$
and divided by a product of elements of $H$.
Counting the number of 
times this term occurs,
 the 
resulting expression
is finally equal to
$$\frac{\binom{n_x-p_x-1}{k_x-p_x-1}y^{n_x-k_x}}{D_x}$$ 
where $D_x$ is a product of $p_x$ linear terms.

Let us prove the last assertion of the proposition. If all the $p_i$s are 0, it means that there is some gluing of $X$ such that for all vertices $x$ of $X$, the resulting cohomology classes $x_1,\ldots,x_{k_x}$ are identical. It means that all edges in $X$ correspond to a null cohomology class, and hence $X$ is a tree.
\end{proof}

\subsection{Application to the integral of torus knots}

We recall that in \cite{lift}, the invariant $\Zwr$ is defined such that $\hair \Zwr(K)=\langle\Omega,\Omega\rangle Z^{\fourwheel}(K)$ for any knot $K$.

Our aim would be to compute $\Zwr(K_{p,q})$ but we have information only on $\hair Z^{rat\fourwheel}(K_{p,q})$. Hence we will compute it only up to the kernel of the $\hair$ map.
We have shown that $\langle\Omega,\Omega\rangle Z^{\fourwheel}(K)=s(X_{p,q})$, hence we are interested in finding a rational expression for $s(X_{p,q})$. Using the preceding section and rationality results of \cite{rat}, we will prove the following theorem:

\begin{thm}\label{theoreme}
Let $X_{p,q}$ be the gluing graph of proposition \ref{Xpq}. It is a gluing graph decorated by the parameters $p,q$ and $pq$ associated respectively to the series $f(px),f(qx)$ and $f(pqx)$. 

We construct from $X_{p,q}$ a rational singular diagram $Y^{\rm rat}_{p,q}\in\boB^{\rm rat}_s$ in the following way.

We consider a tree component $T$ of $X_{p,q}$ and at each vertex $x$, we glue a circle. As the diagram is a tree, the ordering has no importance, and the cohomology class of the circle is well defined: we note it $h_x$. If the vertex $x$ has valence $k_x$ and is colored by $n_x\in\{p,q,pq\}$, we decorate the ``bubble-tree'' diagram with the series $$\prod_x \frac{1}{4}\left(\frac{e^{n_x h_x}+1}{e^{n_x h_x}-1}\right)^{(k_x-1)}.$$
Therefore, we have constructed from $T$ an element $T^{\rm rat}\in\boB^{\rm rat}_s$.
We set $Y^{\rm rat}_{p,q}=\sum\limits_{T\text{ in }X_{p,q}} T^{\rm rat}$.

Then we have $\hair z^{\rm rat\fourwheel}(K_{p,q})=\hair Y^{\rm rat}_{p,q}$.
\end{thm}

In figure \ref{develop}, we give the beginning of the expansion of the unwheeled Kontsevich integral of torus knots obtained thanks to the gluing graph $Y^{\rm rat}_{p,q}$. The parameter $x$ corresponds to the cohomology class of the circle on which the series is located.

\begin{figure}[htbp]
\quad\input{develop.pstex_t}
\caption{Expansion of the unwheeled Kontsevich integral of torus knots}
\label{develop}
\end{figure}

\begin{proof}

Our strategy is to apply the results of the preceding section to $X_{p,q}$: we recall that $X_{p,q}$ is a gluing graph whose power series attached to vertices are either $f(px), f(qx)$ or $f(pqx)$, where $f(x)$ is the series $\frac{1}{2}\log\frac{\sinh(x/2)}{x/2}$. We compute that the derivative of $f$ is  $\frac{1}{4}\coth\frac{x}{2}-\frac{1}{2x}$.

Thanks to proposition \ref{subst}, we can conclude that $s(X_{p,q})$ is obtained by decorating graphs by expressions such as $g(\exp(h))/D$ where $g\in\Q(t)$, $h\in H$ and $D$ is a homogeneous polynomial in $H$.

Hence it is natural to define the following space:

Let $H$ be a free abelian group of finite rank. We note $G(H)$ the image of $S^{-1}\Q[\exp(H)]\otimes (H\setminus\{0\})^{-1}\Q[H]$ in $(H\setminus\{0\})^{-1}\Q[[H]]$. This is again a $\boC$-module and hence, we can define a space $\boB'=\boD(G)$. 

\begin{lem}
The map $S^{-1}\Q[\exp(H)]\otimes (H\setminus\{0\})^{-1}\Q[H]\to G(H)$ is an isomorphism of $\boC$-modules. 
We define the degree of $\frac{P}{Q}\in (H\setminus\{0\})^{-1}\Q[H]$ by $\deg(\frac{P}{Q})=\deg(P)-\deg(Q)$. This degree extends to $G(H)$ and to $\boB'$.
\end{lem}
\begin{proof}
This statement is purely algebraic, it comes from the fact that polynomials in elements of $H$ and exponentials of elements of $H$ are algebraically independant as power series in $H$.
\end{proof}

But we know from \cite{rat} that there is a diagram $z^{\rm rat \fourwheel}(K_{p,q})$ in $\boB^{\rm rat}$ such that $\hair z^{\rm rat \fourwheel}(K_{p,q}) = s(X_{p,q})$. We interpret this result saying that the series $\hair Z^{rat\fourwheel}(K_{p,q})$ lies in the degree 0 part of $\boB'$.

But in the process of substitution of a connected gluing graph, we add denominators unless all edges of the graph are null-homologous (see proposition \ref{subst}). In the same way, all terms containing fractions in $f'(px), f'(qx)$ and $f'(pqx)$ will have a negative degree.


Keeping diagrams in $X_{p,q}$ with degree 0 and substituting them, we obtain exactly the series $\hair Y^{\rm rat}_{p,q}$ described in the theorem. 
As this series made of exactly all gluing graphs producing degree 0 diagrams in $\boB'$ by substitution, we conclude that $\hair Y^{\rm rat}_{p,q}=s(X_{p,q})=\hair z^{\rm rat\fourwheel}(K_{p,q})$.
This states the theorem.
\end{proof}

There are two natural questions for which we have no answer at the moment:
\begin{quest}
\begin{itemize}
\item
How can we apply a weight system to singular diagrams in order to obtain formulas for the Jones functions of torus knots from $Y^{\rm rat}_{p,q}$?
\item
The series $ z^{\rm rat\fourwheel}(K_{p,q})$ and $Y^{\rm rat}_{p,q}$ agree up to the kernel of the $\hair$ map. Are they equal?
\end{itemize}
\end{quest}
\section{Branched coverings}

A great interest for rational expression of Kontsevich integral comes from its relation with branched coverings. More precisely, if $K_{p,q}$ is the torus knot of parameters $p$ and $q$, and $r$ is an integer, let us note $\Sigma^r(K_{p,q})$ be the pair formed of the cyclic branched covering of $S^3$ of order $r$ over $K_{p,q}$ and the ramification link.

If $r$ is coprime with $p$ and $q$, the ramification locus is a knot, and the underlying 3-manifold is a rational homology sphere, the Brieskorn manifold $\Sigma(p,q,r)$. 

In \cite{lift}, a map $\lift_r$ is described which intertwines rational invariant of the cyclic branched coverings and rational invariant of the initial knot in the following way:
$$Z^{\rm rat\fourwheel}(\Sigma^r(K))= \exp(\frac{\sigma_r(K)}{16}\ThetaGraph)\lift_r Z^{\rm rat\fourwheel}(K).$$
We now study this map in the case of torus knots and prove the following proposition:

\begin{prop}
Call $\Pi_r$ the operator on $\boB^{\rm rat}_s$ which multiplies any diagram $D$ by $r^{-\chi(D)}$ where $\chi$ is the Euler characteristic.
then we have $$\lift_r Y^{\rm rat}_{p,q}=\Pi_r Y^{\rm rat}_{p,q}$$
\end{prop}

\begin{proof}

The $\lift_r$ map was defined only for diagrams decorated by fractions without poles at $r$-roots of unity. As we extended the decorations to all fractions, the definition of $\lift_r$ makes sense for any diagram.
In the definition of the $\lift_r$ map, we need to express all denominators as polynomials of $t^r$. Then, we look at the numerators as a coloring by monomials, which is the same as a linear combination of 1-cohomology classes of the underlying graph. We keep only the classes divisible by $r$ and divide them, then we put back denominators replacing $t^r$ by $t$. Finally we multiply the result by $r$.

This construction is very easy in our case because we only need to know how the map $\lift_r$ acts on derivatives of the fraction $\frac{t^n+1}{t^n-1}$ where $n=p,q$ or $pq$. 
Let us write $h(t)=\frac{t^n+1}{t^n-1} \in \Q(t)$. The operator $\lift_r$ and the derivation operator $Dg(t) =tg'(t)$ act on the space $\Q(t)$.

We can develop $h$ in formal series: $h(t)=-1-2\sum\limits_{k\ge 1}t^{nk}$.
This expression shows that we have 
$$\lift_r h(t)=-1 -2\sum\limits_{k\ge 1, r|nk}t^{nk/r}=h(t).$$ 
because $n$ is coprime with $r$.

Then, for $i>0$, we also have $D^i h(t)=-2\sum\limits_{k\ge 1}(nk)^it^{nk}$.
We check in the same way the following formula:
$$\lift_r D^i h(t) = -2\sum\limits_{k\ge 1, r|nk}(nk)^it^{nk/r}=r^i D^i h(t).$$
This shows finally that $\lift_r$ acts on a diagram of $Y^{\rm rat}_{p,q}$ by multiplying it by $r$ to the power $1+\sum_i(v_i-1)$ where $v_i$ is the valence of the $i$-th vertex. This expression is the number of vertices of the diagram minus 1, hence the number of loops minus one.

This ends the proof of the proposition.
\end{proof}

Of course, we have not proved that $\lift_r z^{\rm rat\fourwheel}(K_{p,q})=\Pi_r z^{\rm rat\fourwheel}(K_{p,q})$. Although the first diagram we know in the kernel of the $\hair$ map has loop degree 17 (see \cite{ber}), we do not know if the $\hair$ map is injective in degree greater than 3.

\Addresses\recd

\end{document}